\documentclass[10pt]{ifacconf}

\usepackage{natbib}            
\usepackage{graphicx}          
\usepackage{graphicx}         
\usepackage{amsfonts}
\usepackage{amsmath}
\usepackage{psfrag}

\newcommand{\reals}{\mathbf R}

\usepackage{tikz}
\usetikzlibrary{shapes,arrows}
\usepackage{pgfplots}
\pgfplotsset{plot coordinates/math parser=false}

\begin{document}
\begin{frontmatter}

\title{ Identification of Stochastic Wiener Systems using Indirect Inference\thanksref{footnoteinfo}}

\thanks[footnoteinfo]{This work was partially supported by the Swedish Research Council and the Linnaeus Center ACCESS at KTH. The research leading to these results has received funding from The European Research Council under the European Community's  Seventh Framework program (FP7 2007-2013) / ERC Grant Agrement N. 267381}

\author[First]{Bo Wahlberg}
\author[Second]{James Welsh}
\author[Third]{Lennart Ljung}

\address[First]{Department of Automatic Control and ACCESS, School of Electrical Engineering, KTH Royal Institute of Technology, SE-100 44 Stockholm, Sweden.
 (e-mail: bo.wahlberg@ee.kth.se).}
\address[Second]{School of Elect Engineering and Computer Science, The University of Newcastle, Callaghan NSW 2308, Australia.}
\address[Third]{Division of Automatic Control, Link\"oping University, SE-581 83 Link\"oping, Sweden.}

\begin{abstract}                
We study identification of stochastic Wiener dynamic systems using so-called indirect inference. The main idea is to first fit an auxiliary model to the observed data and then in a second step, often by simulation, fit a more structured model to the estimated auxiliary model. This two-step procedure can be used when the direct maximum-likelihood estimate is difficult or intractable  to compute. One such example is the identification of stochastic Wiener systems, i.e.,~linear dynamic systems with process noise where the output is measured using a non-linear sensor with additive measurement noise. It is in principle possible to evaluate the log-likelihood cost function using numerical integration, but the corresponding optimization problem can be quite intricate. This motivates studying consistent, but sub-optimal, identification methods for stochastic Wiener systems. We will consider indirect inference using the best linear approximation as an auxiliary model. We show that the key to obtain a reliable estimate is to use uncertainty weighting when fitting the stochastic Wiener model to the auxiliary model estimate. The main technical contribution of this paper is the corresponding asymptotic variance analysis.
A numerical evaluation is presented based on a first-order finite impulse response system with a cubic non-linearity, for which certain illustrative analytic properties are derived.
\end{abstract}


\end{frontmatter}

\section{Introduction}
\label{sec:1}
The idea of using two-stage identification methods, e.g., indirect inference, is by no means new in system identification.  Typically, a flexible auxiliary model is fitted to data, and then in a second step this  estimated model is used to find a more structured model. A well-known example of such an approach is the indirect Prediction Error Minimization (PEM) method, \cite{Soderstrometal:91}, where it is assumed that the model structure of interest can be embedded in a larger model structure for which the identi\-fication problem is more tractable. In a second step, the structured model is estimated from the larger model using a weighted non-linear least squares method. The indirect PEM method will, under certain assumptions, have the same asymptotic statistical properties as the Maximum-Likelihood (ML) method, but it can be more efficiently calculated. Indirect PEM is a special case of indirect inference, which was introduced in econometrics  in \cite{gourieroux1993indirect}. Their main motivation was identification problems for which the ML method is intractable. They also proposed the use of Monte Carlo simulations to generate the cost to be minimized in the second step.

\begin{figure}[htb]
\centering
\pgfplotsset{width=1\columnwidth,height=0.5\columnwidth,compat=newest,plot coordinates/math parser=false}\tikzstyle{block} = [draw, rectangle,
    minimum height=3em, minimum width=3em,align=center,node distance=2cm,fill = blue!10]
\tikzstyle{sum} = [draw, circle, node distance=2cm,fill = blue!10]
\tikzstyle{input} = [coordinate]
\tikzstyle{output} = [coordinate]
\tikzstyle{corner} = [coordinate]
\tikzstyle{pinstyle} = [pin edge={to-,thin,black}]

\begin{tikzpicture}[auto, node distance=1cm,>=latex',transform shape,scale=0.85]

    \node [input] (input_u) {};
    \node [block,right of=input_u] (G) {$G(q)$};
    \node [sum, right of=G] (sum1) {$+$};
    \node [input,above of=sum1] (input_w) {};
    \node [block,right of=sum1] (f) {$f(\cdot)$};
    \node [sum, right of=f] (sum2) {$+$};
    \node [input, above of=sum2] (input_e) {};
    \node [output,right of=sum2,node distance=1.5cm] (output_y) {};
        
    \draw [draw,->] (input_u) -- node[above]{$u(t)$} (G);
    \draw [draw,->] (G) -- node[above]{} (sum1);
    \draw [draw,->] (input_w) -- node[above,at start]{$v(t)$} (sum1);
    \draw [draw,->] (sum1) -- node[above]{$z(t)$} (f);
    \draw [draw,->] (f) -- (sum2);
    \draw [draw,->] (input_e) -- node[above,at start]{$e(t)$} (sum2);
    \draw [draw,->] (sum2) -- node[above]{$y(t)$} (output_y);
    
\end{tikzpicture}
\caption{Stochastic Wiener System}\label{fig:1}
\end{figure}
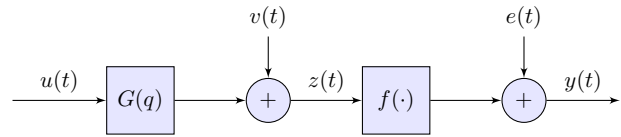

The concept of indirect inference was introduced to the system identification community by \cite{welsh2009continuous} and \cite{larsson2010identification}. The aim of the current paper is to provide further insights into the potential use of indirect inference for the identification of scalar discrete time stochastic Wiener systems, illustrated in Figure \ref{fig:1}, of the form
\begin{align}
z(t)&=G(q)u(t)+v(t),\nonumber\\
y(t)&=f(z(t))+e(t),
 \label{eq:1}
\end{align}
with a stable transfer function $G(q)$ (where $q$ denotes the shift operator), an input signal sequence $\{u(t)\}$, white  stationary process noise $\{v(t)\}$ with zero mean and variance $\sigma_v^2$,  an output signal sequence $\{y(t)\}$, and additive white stationary measurement
noise $\{e(t)\}$ with zero mean and variance $\sigma_e^2$. The input signal $u(t)$ is assumed to be independent of the noises (the open-loop case), and may be a realization of a stochastic process with known probability distribution. To simplify the presentation, we will assume that the noise processes are independent and normal (gaussian) distributed. The more general case with coloured process noise can be handled using a predictor model.

The main challenge is the non-linear function $f(\cdot)$, which means that we have a non-linear stochastic system where the process noise $\{v(t)\}$ propagates through a non-linear device. This typically corresponds to a non-linear sensor.

\subsection{Contribution}
The main objective of this paper is to study the indirect inference method using the best linear approximation as an auxiliary model and compare this to maximum-likelihood  and prediction error minimization methods for identification of a stochastic Wiener system.  The main challenge is to handle the non-linear process noise contribution. It is in principle possible to calculate the likelihood cost, or some approximation of it, using for example particle filters and/or Markov Chain Monte Carlo methods, but even for simple examples the maximum-likelihood approach leads to rather involved computations. A problem for prediction error minimization methods  is that the probability density function of the conditional mean prediction error is quite complicated and involves convolution integrals. This is the reason for using more ad hoc identification techniques. A common approach is to use linear and gaussian approximations. The main technical contributions are:
\begin{itemize}
\item To connect the use of best linear approximation of stochastic Wiener systems  with the method of indirect inference. We present the corresponding variance uncertainty weighting, which for this case includes the input signal. We also derive the asymptotic variance expression for the estimated model parameters.
\item Illustration of the results on a first-order finite impulse response model, for which it is possible to find  analytic expressions. The statistical performance on this example of indirect inference based on best linear approximation is comparable to that that of maximum-likelihood estimation and prediction error minimization. The computations for the indirect inference method are just a fraction of the ones for calculating the maximum-likelihood estimate. The example also shows  that the cost of using a non-linear sensor is increased uncertainty in the estimated model.
\end{itemize}
The statistical theory for identification of stochastic Wiener systems is by no means complete, and  the aim of this paper is to provide insights in some open  important problems.

\subsection{System Identification}
Given measurements of the input and output signals  $\{u(t), y(t)\}$, $t=1,\ldots,N$, the task is to identify a model of the stochastic Wiener system of the form,
 \begin{align}
z(t)&=G(q,\theta)u(t)+v(t),\nonumber\\
y(t)&=f(z(t))+e(t)
 \label{eq:2}
\end{align}
where the model is parameterized by the parameter vector $\theta\in \reals^n$. We assume that the true system can be described by $\theta_o$. The noise processes $\{v(t)\}$ and $\{e(t)\}$ are assumed to be independent normal distributed (gaussian) zero mean white noise.
The corresponding noise variances $\sigma_v^2$ and $\sigma_e^2$ are assumed to be known, but could be added to the parameters to be estimated. This can, however, cause identifiability problems.  We will  study the case when the non-linear function $f(\cdot)$ is known. It is possible to extend our results to the case when the function  $f(\cdot)$ also is estimated, which can result in an identifiability problem. The reason for these simplifications is to focus on the stochastic part due to the process noise $v(t)$.

Identification of Wiener systems is a well studied topic, see for example \cite{149500,wigren1993recursive,bai2003frequency,zhu2002estimation,enqvist2005linear,Pillonetto20132704,Giri:10} and the references therein. It forms the basis for the identification of more general non-linear block diagram based models. However, many algorithms assume no process noise, which leads to a non-linear least squares problem minimizing the output error, i.e.,~the difference between the measured and the simulated outputs.  The maximum-likelihood  method for stochastic Wiener systems was introduced in \cite{914162} and analysed in more detail in \cite{hagenblad2008maximum}. The expectation-maximization  algorithm for maximum-likelihood identification and the use of the particle filter have been studied in \cite{wills2013identification,wills2010wiener}. As recently pointed out in \cite{wahlbergcdc2014}, the stochastic Wiener system identification problem can be viewed as a non-linear errors-in-variables problem, with the well-known bias problem due to input noise. This paper also proposes a prediction error minimization  framework for identification of stochastic Wiener systems. The idea is to use the conditional mean predictor and  notice that the variance of the prediction errors may be  highly dependent on the input signal $u(t)$. Hence, variance uncertainty weighting is most important in order to obtain reliable estimates.

\subsection{Outline}
 ML and PEM methods for stochastic Wiener system identification are reviewed in Sec.~\ref{sec:2}.
In Sec.~\ref{sec:3},  we study how to use the indirect inference approach for the identification of stochastic Wiener models by using the BLA as the auxiliary model. A first order FIR model example is outlined in Sec.~\ref{sec:5}, and the corresponding
numerical study is presented in Sec.~\ref{sec:6}. The paper is concluded in Sec.~\ref{sec:7}.

\section{ML and PEM}
\label{sec:2}
As shown in \cite{914162}, the negative log-likelihood function, given data and the normal distributed noise model (\ref{eq:2}), equals
$$
l(\theta)=-\sum_{t=1}^N \log \int_{-\infty}^\infty e^{-{\cal E} (t,\theta, z)} dz,
$$
where
$$
{\cal E}(t,\theta,z)=\frac{[y(t)-f(z)]^2}{2\sigma_e^2}+\frac{[z-G(q,\theta)u(t)]^2}{2\sigma_v^2}.
$$
The ML estimate of $\theta $ is obtained by minimizing  $l(\theta)$.
There are at least two challenges with the  ML method for stochastic Wiener systems. First, to evaluate the negative log-likelihood cost at a certain parameter value $\theta$ we have to calculate $N$ integrals. This can be done rather efficiently using numerical integration and parallel computations. A difficulty is the integrand, where typically
$$
e^{-{\cal E} (t,\theta, x)}\approx \left\{\begin{array}{cl} 1, & x\; \mbox {small,}\\  0, & x\; \mbox {otherwise.}
\end{array}\right .
$$
The example to be considered in Sec.~\ref{sec:5} corresponds to ${\cal E} (x)\sim x^6$, which means that the integrand decreases very rapidly to zero.

The PEM approach avoids the exponential function integration issue by using a weighted least squares cost function, see
\cite{wahlbergcdc2014}. The conditional mean predictor
 of $y(t)$ for given $u(t)$ and $\theta$ is
\begin{align}
\label{eq:LL1}
\hat{y}(t,\theta)& ={\mathrm E}_v\{f (G(q,\theta)u(t)+v(t))\}.
\end{align}
Notice that the prediction error variance depends on the input signal $u(t)$.
The optimally weighted quadratic PEM cost-function, see \cite{wahlbergcdc2014}, is
\begin{align}
V_{N}(\theta)&=\frac{1}{N}\sum_{t=1}^N\frac{\epsilon^2(t,\theta)}{{\mathrm E}\{\epsilon^2(t,\hat{\theta}_I)\}},
\label{eq:wpem}
\end{align}
with prediction error
$
\epsilon(t,\theta)=y(t)-\hat{y}(t,\theta).
$
The variance weighting is calculated  at a consistent initial estimate $\hat{\theta}_I$, e.g.,~the  PEM estimate without weighting.
The use of weighting is important to obtain reliable estimates and depends here on the input signal $u(t)$.
The PEM estimate based on (\ref{eq:wpem}) is not asymptotically efficient for non-linear functions, since we use a weighted quadratic cost-function. However, by using a cost-function based on the probability density function of $\epsilon(t,\theta)$, we obtain an asymptotically efficient PEM estimate, see \cite{ljung1998system}.  The computations will then be similar to the ML case, involving multiple integral calculations with exponential functions.

\section{Indirect Inference Using Best Linear Approximation}
\label{sec:3}
\subsection{Best Linear Approximation of Stochastic Wiener Systems}

Let us illustrate the concept of indirect inference applied  to a stochastic Wiener system with a known non-linear function. A common first approach is to fit a linear model to the observed data. It is well known that if the input signal is normal (gaussian) distributed, then the Best Linear Approximation (BLA), see \cite{Ljung:01} or \cite[Chapter 13]{Giri:10},  is a scaled version of the linear dynamics transfer function $G(q)$ of the Wiener system. It is perhaps less well known that the same result holds for {\em stochastic} Wiener systems with gaussian process noise $v(t)$. This extension follows more or less from Bussgang's theorem, \cite{Bussgang52}: If $z(t)$ is a normal distributed stationary stochastic process with zero mean and  if the non-linear  transformed process $y(t)=f(z(t))$ has zero mean, then
$$
{\mathrm
  E}\{y(t)z(t-\tau)\}=b_0 {\mathrm E}\{z(t)z(t-\tau)\},\quad b_0={\mathrm E}\{f'(z(t))\}.
$$
Adding independent gaussian process noise to the input signal contribution,  $z(t)=G(q)u(t)+v(t)$, makes $z(t)$ still normal distributed and Bussgang's theorem holds.  Now we are only interested in computing the partial correlations using the relations
\begin{align}
{\mathrm
  E}\{y(t)u(t-\tau)\}&=b_0 {\mathrm E}\{z(t)u(t-\tau)\}\nonumber\\ &= b_0 {\mathrm E}\{[G(q)u(t)]u(t-\tau)\}.
  \label{eq:8}
\end{align}
A recent proof of (\ref{eq:8}) can be found in \cite{Banelli:12}, but it already follows from \cite{Nutall58}.
This result proves that the BLA of a stochastic Wiener system with a normal distributed input signal again equals
$$
G_{BLA}(q)= b_0 G(q)
$$
i.e.,~a scaled version of the linear transfer function. Notice that $b_0$ now also depends on the statistics of the process noise $v$ as well as of the input $u$.
This result can be generalized to separable stochastic processes, \cite{Nutall58,enqvist2005linear}.

The BLA of a stochastic  Wiener system for more general input signal sequences can be obtained by simulations or by analytic calculations (if the distribution of the input signal is known).

\subsection{An Optimally Weighted  Indirect Inference Algorithm}
We will now describe how to apply indirect inference to the stochastic Wiener system identification problem using a two-step procedure based on BLA. From the given data   $\{u(t), y(t)\}$, $t=1,\ldots,N$, form the BLA cost-function
\begin{equation}
Q_N( \beta)=\frac{1}{N}\sum_{t=1}^N[y(t)-G_{lin}(q,\beta)u(t)]^2.
\label{eq:q}
\end{equation}
Here we have used the notation $G_{lin}(q,\beta)$ to stress that this is a  linear transfer function  parameterized by $\beta\in \reals^m$. It is in general of higher order compared to $G(q,\theta)$ in (\ref{eq:2}), that is $m$ is typically larger than $n$, the dimension of $\theta$. We have used an Output Error PEM cost-function (\ref{eq:q}), but it is also possible to use a more general PEM model structure, as described in \cite{Ljung:01}.

{\bf Step 1:} Identify the BLA within the model structure $G_{lin}(q,\beta)$, $\beta\in \reals^m$, of the  system from given data using the cost-function (\ref{eq:q})
$$
\hat{\beta}_N= \mbox{arg}\min_{\beta} \{Q_N( \beta)\},
$$
i.e.,~the BLA estimate equals $\hat{G}_{BLA}(q)=G_{lin}(q,\hat{\beta}_N)$.

The next question is to figure out the functional relation  $\beta(\theta)$, $\reals^{m}\mapsto \reals^{n}$ that describes how the auxiliary estimate $\hat{\beta}_N$ asymptotically depends on the underlying structured model of interest, as a function of the  model parameter vector $\theta$. Now assume that $Q_N( \beta)$ converges (almost surely) as the number of data $N\to \infty$ to $Q( \beta,\theta_o)$, which depends on the true system parameter vector $\theta_o$.  Here
\begin{align*}
Q( \beta,\theta)&= {\mathrm E}\{[y(t)-G_{lin}(q,\beta)u(t)]^2\}\nonumber\\
&={\mathrm E}_{v,u}\{[f(G(q,\theta)u(t)+v(t))-G_{lin}(q,\beta)u(t)]^2\}.
\end{align*}
The expectation is with respect to both the input signal and the process noise. It is also be possible to only use expectation with respect to the process noise for a given input sequence.
Let
$$
\beta(\theta)= \mbox{arg}\min_{\beta} \{Q(\beta,\theta)\},
$$
and use the notation
$$
\beta_o=\beta(\theta_o),
$$
for the corresponding true parameter vector.
This leads to:

{\bf Step 2: (Analytic)} Estimate the structured model parameter vector $\theta$ by solving
\begin{equation}
\hat{\theta}_{N}=\mbox{arg}\min_{\theta}[{\beta}(\theta)-\hat{\beta}_N]^TW [{\beta}(\theta)-\hat{\beta}_N],
\label{eq:thetaanalytic}
\end{equation}
where $W$ is a positive definite weighting matrix to be specified.

For certain examples of  non-linear functions and distributions it may not possible to analytically find the function $\beta(\theta)$. One can then resort to  Monte Carlo simulations. Let
\begin{align*}
&Q_{N,S}( \beta,\theta)= \nonumber\\
&\frac{1}{S}\sum_{s=1}^S\frac{1}{N}\sum_{t=1}^N[f(G(q,\theta)u(t)+v_s(t))-G_{lin}(q,\beta)u(t)]^2\}.
\end{align*}
Here $\{v_s(t)\}$, $t=1,\ldots , N$, $s=1. \ldots, S$, is a generated noise realization of the process noise $v(t)$ and $S$ is the total number of realizations used in the Monte Carlo Simulation. Let
$$
\hat{\beta}_{N,S}(\theta)= \mbox{arg}\min_{\beta}\{ Q_{N,S}( \beta,\theta)\}.
$$

{\bf Step 2: (Simulated)} Estimate the structured model parameter vector $\theta$ by solving
\begin{equation}\label{eq:thetasim}
\hat{\theta}_{N,S}=\mbox{arg}\min_{\theta}[\hat{\beta}_{N,S}(\theta)-\hat{\beta}_N]^TW [\hat{\beta}_{N,S}(\theta)-\hat{\beta}_N],
\end{equation}
where $W$ is a  positive definite weighting matrix to be specified.

The  {\em optimal weighting matrix}, $W$, equals the inverse of the covariance matrix of the auxiliary estimate $\hat{\beta}_N$,  compare Expression (9.11) for the asymptotic covariance matrix in \cite{ljung1998system}. In our setting this translates to
\begin{align}
W&=[J^{-1}_oI_oJ_o^{-1}]^{-1},\label{eq:WW}\\\
I_o&= \lim_{N\to \infty}\text{Cov}\{\sqrt{N} \frac{\partial Q_N}{\partial \beta}(\beta_o,\theta_o)\},\;J_o= \frac{\partial^2 Q}{\partial \beta^2}(\beta_o,\theta_o).\nonumber
\end{align}
In practise, one has to use a consistent estimate of $W$.
\subsection{Performance Analysis}

Next, we determine the asymptotic properties of the final estimate $\hat{\theta}_{N}$. This can also be done using Taylor series expansion arguments as described in Complement C4.4  in \cite{Soderstrom&Stoica:89}. The properties of the function ${\beta}(\theta)$, $\reals^{m}\mapsto \reals^{n}$, are very important in order to obtain a consistent estimate. It should be possible to invert this function, $\theta=\alpha({\beta}(\theta))$, i.e.,~$\alpha(\cdot)$ is a left inverse of $\beta(\cdot)$.
Define the Jacobian matrix
\begin{align}
G=\frac{\partial {\beta}}{\partial\theta}(\theta_o)\in \reals^{m\times n},
\label{eq:g}
\end{align}
and recall that we use the weighting matrix
\begin{equation}
W=[\lim_{N\to \infty} \text{Cov}\{\sqrt{N}[\hat{\beta}_N-\beta_o]\}]^{-1}.
\label{eq:w}
\end{equation}
{\em Variance Expression:} The asymptotic covariance matrix of $\hat{\theta}_N$ defined by (\ref{eq:thetaanalytic}) equals
\begin{equation}
\lim_{N\to \infty} \text{Cov}\{\sqrt{N} [\hat{\theta}_N-\theta_o]\}= [G^TWG]^{-1},
\label{eq:main}
\end{equation}
where $G$ is defined by (\ref{eq:g}) and $W$ by (\ref{eq:w}).

To prove this result, we will make a Taylor series expansion of the cost-function at the minimizing $\theta$:
$$
V_w(\theta)=[{\beta}(\theta)-\hat{\beta}_N]^TW [{\beta}(\theta)-\hat{\beta}_N]\;\Rightarrow \; V'_w(\hat{\theta}_N)=0,
$$
which gives  $[\hat{\theta}_N-\theta_o]\approx  -[V''_w({\theta}_o)]^{-1}V'_w(\theta_o).$
The Hessian $V''_w({\theta}_o)$ has to be invertible (positive definite for unique local minimum) for this to hold. As shown  in Complement 4.4 in  \cite{Soderstrom&Stoica:89}
\begin{align*}
V'_w({\theta}_o)&=2G^TW[\hat{\beta}_N-\beta_o]+{\cal O}(1/{N})\\
V''_w({\theta}_o)&=2G^TWG+{\cal O}(1/\sqrt{N}).
\end{align*}
Hence
$$
[\hat{\theta}_N-\theta_o]\approx  -[G^TWG]^{-1}G^TW[\hat{\beta}_N-\beta_o].
$$
Since $\lim_{N\to \infty} \text{Cov}\{\sqrt{N}[\hat{\beta}_N-\beta_o]\}=W^{-1}$, this proves~(\ref{eq:main}).

When using the Monte Carlo simulation approach in Step~2, (\ref{eq:thetasim}), the asymptotic covariance matrix of $\hat{\theta}_{N,S}$ should be amplified by the factor
$(1+{1}/{S})$, due to the extra uncertainty from the simulations, see \cite{heggland2004estimating}.

\subsection{Comments}

The suggestion to use the BLA  as an auxiliary model is rather ad hoc, but a very common choice in recent methods for identification of non-linear systems, \cite{schoukens2005identification,pintelon2012system,schoukens2014identification,schoukens2012parametric,Sjoberg12}.   \cite{Ljung:01} contains an overview of the role of BLA in system identification.
Another option would be to use the biased estimate from minimizing
$$
l(\beta)=\sum_{t=1}^N [y(t)-f(G(q, \beta)u(t)]^2
$$
and then use Step 2 to remove the bias of this estimate.
A challenge is to find even more efficient auxiliary models. The key property is that ${\beta}(\theta)$ should have a left inverse (identifiability) and
\begin{align}
G=\frac{\partial {\beta}}{\partial\theta}(\theta_o)\in \reals^{m\times n}
\end{align}
should be     ``large'' (sensitive) and at the same time $\beta$ should be easy to estimate from the given data.

\subsection{Indirect Inference}
\label{sec:4}

As mentioned in the introduction the indirect inference approach was developed in \cite{gourieroux1993indirect} as a way to find
consistent model estimates even when the ML method is intractable. Our description of the indirect inference approach applied to BLA in the previous section is mainly based on \cite{heggland2004estimating}.  The theory of indirect inference is more general, and the key step is to choose the auxiliary model parameterized by $\beta$ and the data-driven cost-function
$Q_N(\beta)$. A  convergence analysis of the general indirect inference method  presented in  \cite{heggland2004estimating} is based on using the concept of finite dimension auxiliary statistics.

If  $Q_N(\beta)$ is based on a  {\em sufficient statistics} the indirect inference method is efficient (except for the factor $(1+1/S))$. This is not the case for stochastic Wiener systems, for which  no finite dimensional sufficient statistics exists. For the BLA approach we are only using a second order statistics.

\section{Illustrating Example}
\label{sec:5}

We will use the following simple example to try to further understand the properties of the stochastic Wiener system identification methods described in the previous sections,
\begin{align}
z(t)&=\theta u(t)+u(t-1)+v(t)\nonumber\\
y(t)&=[z(t)]^3+e(t).
\label{eq:simsys}
\end{align}
The same example was used in \cite{larsson2010identification,wahlbergcdc2014}. We will study two different input white noise distributions, a normal (gaussian) and an uniform distribution, both with variance $\sigma_u^2$. The noises are assumed to be white zero mean normal distributed with variances $\sigma_e^2$ and $\sigma_v^2$, respectively.

The maximum-likelihood cost-function
$$
l(\theta)=-\sum_{t=1}^N \log {\mathrm E}_{\bar{v}}\{\frac{\sigma_v}{\sqrt{2\pi}} e^{-\frac{1}{2\sigma_e^2}[y(t)-(\theta u(t)+u(t-1)+\sigma_v\bar{v})^3]^2}\}
$$
is calculated using a Gauss Hermite approximation of order $1000$. The reason for this high order is the $e^{-x^6}$ tends to zero very quickly and the integral is rather difficult to approximate. We  have noticed that increased process noise variance $\sigma_v^2$ makes it even more challenging.

The conditional mean predictor equals
\begin{align}
\hat{y}(t,\theta)&=\theta^3 u^3(t)+u^3(t-1)+3(\theta u^2(t)u(t-1)
\nonumber\\ &+\theta u(t)u^2(t-1)+\theta u(t)\sigma_v^2 +u(t-1)\sigma_v^2),
\end{align}
and is used in the PEM method (\ref{eq:wpem}).

In order to evaluate the indirect inference method we will study the BLA. Assume first that the input signal is normal distributed with zero mean and variance $\sigma_u^2$. The BLA of model order one, i.e.,~based on minimizing
$$
{\mathrm E} \{(y(t)-\beta_1 u(t)+\beta_2u(t-1))^2\}
$$
is
\begin{equation}
\left(\begin{array}{c} \beta_1(\theta)\\ \beta_2(\theta)\end{array}\right)= [3\sigma_u^2\theta^2+3(\sigma_u^2+\sigma_v^2)]\left(\begin{array}{ccc} \theta \\
 1
 \end{array}\right),
 \label{eq:g1}
\end{equation}
for which
$$
\frac{\beta_1}{\beta_2}=\theta,
$$
since the BLA for gaussian input signal is just a scaled version of $G(q,\theta)$.  This also gives that $\beta(\alpha)$ has a left inverse and we have identifiability of $\theta$ from $\beta$.

The corresponding BLA when the input signal is uniformly distributed with zero mean variance $\sigma_u^2$ is
\begin{equation}
\left(\begin{array}{c} \beta_1(\theta)\\ \beta_2(\theta)\end{array}\right)=\left(\begin{array}{ccc} \frac{9}{5}\sigma_u^2\theta^3+3(\sigma_u^2+\sigma_v^2)\theta \\
 3\sigma_u^2\theta^2+3(\frac{3}{5}\sigma_u^2+\sigma_v^2)
 \end{array}\right).
\label{eq:g2}
\end{equation}
For the indirect inference method we will use
$$
Q_N(\beta)=\frac{1}{N}\sum_{t=1}^N(y(t)-\beta_1 u(t)+\beta_2u(t-1))^2
$$
and hence
$$
\frac{\partial Q_N}{\partial \beta}=\frac{2}{N}\sum_{t=1}^N \left(\begin{array}{c}  (y(t)-\beta_1 u(t)+\beta_2u(t-1))u(t)\\
(y(t)-\beta_1 u(t)+\beta_2u(t-1))u(t-1)
\end{array}\right)
$$ and we need to calculate the weighting matrix $W$ from (\ref{eq:WW}). Here $J_o=\sigma_u^2I $ and the tedious work is to calculate $I_o$.
 It can, however, be estimated as
\begin{align*}
\hat{I}_o&=\frac{1}{N}\sum_{t=1}^N\left[[y(t)-\hat{\beta}_1 u(t)-\hat{\beta}_2u(t-1)]^2\right.\times\nonumber\\ &\left. \left(\begin{array}{cc}  u^2(t) & u(t)u(t-1)\\
u(t))u(t-1) & u^2(t-1)
\end{array}\right)\right]
\end{align*}
\vskip0.5\baselineskip
{\bf Comment:} The key to get the order one indirect inference method to work is to use the weighting $W=\hat{I}_o^{-1}$.  The paper \cite{larsson2010identification} does not use weighting, which explains their non-intuitive simulation result that it is better to use a zero order BLA model than a first order one (which should contain more information).

We will finally study the zero order model case, for which no weighting is needed. Here
$$
Q_N(\beta)=\frac{1}{N}\sum_{t=1}^N(y(t)-\beta_1 u(t))^2
$$
with the BLA given by $\beta_1(\theta)$ in (\ref{eq:g1}) and (\ref{eq:g2}).
An alternative cost-function is
$$
Q_N(\beta)=\frac{1}{N}\sum_{t=1}^N(y(t)-\beta_2 u(t-2))^2,
$$
with the BLA given by $\beta_2(\theta)$ from (\ref{eq:g1}) or (\ref{eq:g2}).
A problem here is that $\beta_2(\theta)$ is quadratic in $\theta$ and we do not have an unique solution with respect to $\theta$.

\section{ Simulation Result}
\label{sec:6}
We will use the following numerical values for the system parameters in (\ref{eq:simsys}):
$$\theta_o=0.5,\quad\sigma_v^2=0.2,\quad \sigma_e^2=0.1.$$ We will use the analytic Step~2, (\ref{eq:thetaanalytic}), for the indirect inference based on BLA, and will evaluate the following methods for $N=1000$ observations:

\begin{align}
\mbox{Method 1:}& \quad \mbox{ML}\nonumber\\
\mbox{Method 2:}& \quad \mbox{PEM with optimal weighting}\nonumber\\
\mbox{Method 3:}& \quad \mbox{Zero order indirect inference}\nonumber\\
\mbox{Method 4:}& \quad \mbox{First order indirect inference, no weighting }\nonumber\\
\mbox{Method 5:}& \quad \mbox{First order indirect inference, with weighting\nonumber}
\end{align}

The following table summarizes the simulation results for the {\em normal distributed input sequence} with zero mean and variance $\sigma^2_u=1/3$
over $1000$ noise  and input realizations:

\begin{align}
\mbox{Method 1:}& \quad \mbox{mean: $0.5025$ \quad std: $0.0303$ }\nonumber\\
\mbox{Method 2:}& \quad \mbox{mean: $0.5014$ \quad std: $0.0349$ }\nonumber\\
\mbox{Method 3:}& \quad \mbox{mean: $0.4983$ \quad std: $0.0446$ }\nonumber\\
\mbox{Method 4:}& \quad \mbox{mean: $0.4977$ \quad std: $0.0554$ }\nonumber\\
\mbox{Method 5:}& \quad \mbox{mean: $0.4982$ \quad std: $0.0418$}\nonumber
\end{align}

The next table summarizes the simulation results for the {\em uniform distributed input sequence } with zero mean and variance $\sigma^2_u=1/3 $ over $1000$ noise  and input realizations

\begin{align}
\mbox{Method 1:}& \quad \mbox{mean: $0.4994$ \quad std: $0.0325$ }\nonumber\\
\mbox{Method 2:}& \quad \mbox{mean: $0.4984$ \quad std: $0.0346$ }\nonumber\\
\mbox{Method 3:}& \quad \mbox{mean: $0.4988$ \quad std: $0.0454$ }\nonumber\\
\mbox{Method 4:}& \quad \mbox{mean: $0.4984$ \quad std: $0.0458$ }\nonumber\\
\mbox{Method 5:}& \quad \mbox{mean: $0.4987$ \quad std: $0.0377$}\nonumber
\end{align}

The conclusion from the numerical study is that the ML, as expected, gives the best performance. However, PEM and indirect inference with optimal weighting give also good results. The computational times for these methods are only a fraction of that of ML.  The simulation study also shows that a zero order model can give better results than using a first order model and no weighting, c.f.,~\cite{larsson2010identification}.

It should be noted that the asymptotic variance when $f(x)=x$ equals
$
(\sigma_v^2+\sigma_e^2)/(\sigma_u^2 N),
$
with a standard deviation equal to $0.0173$ for the example above. Hence, the non-linear function $f(x)=x^3$ gives a more difficult estimation problem than the standard case $f(x)=x$.

\section{Conclusion}
\label{sec:7}

In this paper, we have utilised the indirect inference method based on the best linear approximation to identify stochastic Wiener systems. The results show that to obtain a reliable estimate, it is important to use an optimal weighting when estimating the structured model from the auxiliary model. The weighting here is the inverse of the covariance matrix of the BLA parameter estimate. We have analyzed the statistical properties of the corresponding indirect inference estimate using results from system identification. The methods have been evaluated using a first order FIR system with a cubic non-linearity. The simulation results demonstrate that the proposed indirect inference BLA approach performs quite well compared to the ML and weighted PEM methods. A major advantage is that the indirect inference algorithms are very computationally fast and direct to implement.

There are many open questions when it comes to identification of stochastic non-linear dynamic systems. For example, performance results to guide the choice of sensors.


\bibliography{refs,referenslista}

\end{document}